\documentclass[11pt]{article} 
\usepackage[utf8]{inputenc}
\usepackage{caption,graphicx,amsmath,amsthm,amssymb, authblk,xcolor,enumitem,subfigure,float,verbatim,hyperref,pgfplots,bibentry,mathtools,soul,isomath, cancel,multicol,wrapfig,subcaption}
\usepackage{mathrsfs}
\usepackage{indentfirst}
\usepackage{nccmath}
\graphicspath{ {images/} }
\usepackage{multicol}
\usepackage{graphicx,amsmath,amssymb,amsthm,marvosym,tikz,fontenc}
\usetikzlibrary{backgrounds}
\usetikzlibrary{patterns,shapes.misc, positioning,fit,calc}
\usepackage{tikz}
\usetikzlibrary{graphs,graphs.standard}
\usetikzlibrary{decorations.pathmorphing}
\usepackage{pgfplots}
\pgfplotsset{width=10cm,compat=1.9,tick scale binop=\times}
\usepackage{indentfirst}
\usepackage[a4paper,bindingoffset=0.2in,%
left=0.8in,right=0.8in,top=1.2in,bottom=1.2in,%
footskip=.25in]{geometry}
\usepackage{relsize}
\usepackage[english]{babel}
\allowdisplaybreaks
\theoremstyle{plain}
\newtheorem{theorem}{Theorem}[section]
\newtheorem{lemma}[theorem]{Lemma}

\newtheorem{corollary}[theorem]{Corollary}
\newtheorem{definition}[theorem]{Definition}
\newtheorem{remark}[theorem]{Remark}
\newtheorem{example}[theorem]{Example}

\newlength{\defbaselineskip}
\newcommand{\setlinespacing}[1]%
{\setlength{\baselineskip}{#1 \defbaselineskip}}

\date{}
\begin{document}
\title{On the Estrada Index and Spectral Properties of $k$-Uniform Hypergraphs}
\author{Liya Jess Kurian$^1$\footnote{liyajess@gmail.com},  Chithra A. V$^1$\footnote{chithra@nitc.ac.in}
 \\ \small 
 1 Department of Mathematics, National Institute of Technology Calicut,\\\small
 Calicut-673 601, Kerala, India\\ \small
 $^*$liyajess@gmail.com, $^\dagger$Corresponding author email: chithra@nitc.ac.in}
\maketitle
\thispagestyle{empty}
\begin{abstract}
Various properties of a hypergraph can be explored through its spectrum. In this paper,we estimate the bound for the sum of $t,\,t\in[1,n]$, largest eigenvalues of a $k$-uniform hypergraph of order $n$. Also, we characterize the $k$-uniform hypergraph with two distinct eigenvalues. We establish bounds for the Estrada index of $k$-uniform hypergraphs based on their number of edges, order and energy. In addition, we obtain the second largest Estrada index among the unicyclic hypergraphs.
\end{abstract}

\section{Introduction}
Spectral graph theory provides powerful tools for analyzing the structure of graphs using the eigenvalues and eigenvectors of the matrices associated with them. It is natural to extend the theory of spectral graphs theory from graphs to hypergraphs. In 2008, Lim \cite{lim2008} proposed the study of hypergraph spectra by means of tensors. The spectral analysis of hypergraphs through tensors is computationally challenging. Recently, researchers have shifted their focus to matrix-based spectral methods for hypergraphs\cite{Banerjee2021, Kumar2017, Lin2017}. 

Throughout this article, we consider simple and finite hypergraphs. More precisely, a hypergraph $G^*=(V(G^*),E(G^*))$ of order $n$ with vertex set $V(G^*)=\{v_1,v_2,v_3,\ldots,v_n\}$ and hyperedge set (edge set) $E(G^*)=\{e_1,e_2,e_3,\ldots,e_m\}$ which is the collection of hyperedges (edges) such that each hyperedge $e_i,\,|e_i| \geq 2$ is a non-empty subset of the vertex set $V(G^*)$. The degree $d(v)$ of a vertex $v\in V(G^*)$ is defined as the number of edges that contain the vertex $v$. A hypergraph $G^*$ is said to be a $k$-uniform hypergraph if each of its edges contains exactly $k$  vertices \textnormal{\cite{Cooper2012,Kumar2017}}, and $r$-regular if $d(v_i)=r$ for all $v_i\in V(G^*)$. A hypergraph is said to be a $(k,r)$-regular hypergraph if it is both  $k$-uniform  and $r$-regular. Specifically, a complete $k$-uniform hypergraph $(K_n^k)$ of order $n$ is a $\left(k,\binom{n-1}{k-1}\right)$-regular hypergraph. Its complement $\overline{K}_n^k$ is the edgeless hypergraph in which no two vertices are connected. The adjacency matrix $A(G^*)=(a_{ij})$ of $G^*$ \cite{Lin2017} is an $ n\times n $  matrix whose rows and columns are indexed by the vertices of $G^*$ and for all $ v_i,v_j \in V, $
 \begin{equation*}
a_{ij} =\left\{
   \begin{array}{ll}
      \mid \{e_k  \in E(G^*): \{v_i,v_j\} \subset  e_k\}\mid    &  \mbox{, } v_i \neq v_j  \\
       0  & \mbox{, } v_i=v_j
   \end{array}.
   \right.
\end{equation*}

For a hypergraph of order \( n \), let \(\lambda_i, \,i\in[1,n]\) be the eigenvalues of \( A(G^*) \), where \( \lambda_1 \geq \lambda_2 \geq \cdots \geq \lambda_n \). In this paper, we estimate an upper bound for sum of $t$- largest eigenvalues ($\sum_{i=1}^t\lambda_i,t\in[1,n]$) of a hypergraph $G^*$. 
Also, this paper examines $k$-uniform hypergraphs with two distinct eigenvalues in relation to the combinatorial properties of interconnection networks. Various networks can be effectively modeled using hypergraphs. 

 One of the important spectrum-based graph invariants is the Estrada index, which was first introduced by Estrada in 2000\cite{estrada2000}.  It has various applications in chemistry and physics. 
In other words, we can write it as
\begin{equation}\label{eqn1a}
    EE(G^*)=\sum_{t\geq0}\dfrac{M_t(G^*)}{t!}=\sum_{t\geq0}\dfrac{\mathrm{tr}(A(G^*))^t}{t!}=\mathrm{tr}(e^{A(G^*)}),
\end{equation}
 where $M_t(G^*)$ is the $t$-th spectral moment of hypergraph $G^*$ defined as $M_t(G^*)=\sum_{i=1}^n\lambda_i^t$ and $\mathrm{tr}(A(G^*))$ represents the trace of the matrix $A(G^*)$. Many interesting results have been obtained while analyzing hypergraphs that maximize or minimize the Estrada index within a given class of hypergraphs. For example, the 3-uniform hypergraph with the maximum Estrada index has been identified\cite{li2014}, as well as the r-uniform linear hypertrees with the extremal Estrada index\cite{Wand2020}. Moreover, the unicyclic hypergraph that attains the maximum Estrada index has also been determined\cite{feng2023}. In this paper, we obtain the second largest Estrada index among the unicyclic hypergraphs.\\
This paper focuses on the study of the spectral properties of hypergraphs and the Estrada index associated with the hypergraphs. In Section 2, we present the fundamental definitions, notations, and results that will be used throughout this paper.  In Section 3, we estimate an upper bound for the sum of the t-largest eigenvalues of a $k$-uniform hypergraph $G^*$. Also, we characterize the $k$-uniform hypergraphs with two distinct eigenvalues. In Section 4, we determine the bounds for the Estrada index of the $k$-uniform hypergraph in terms of its order, size, and energy. In addition, we provide Nordhaus-Gaddum type bound for the Estrada index of $k$-uniform hypergraphs. Furthermore, we identify the $k$-uniform hypergraph with the second largest Estrada index within the unicyclic 
$k$-uniform hypergraph.
\section{Preliminaries}
This section gives the basic definitions, results, and notation used in the main results.\\

We write $j\in [a,b]$ if $j$ takes all the integer values satisfying the condition $a\leq j\leq b$ and for any set $S$, $\mathcal{P}(S)$ denote the power set of $S$. We denote all one matrix and identity matrix of order $n$ by $J_n$ and $I_n$, respectively. Let \( S_n \) be the set of all real symmetric matrices of order \( n \) and
\(M_n = \{ M=(m_{ij}) \in S_n : 0 \leq m_{ij} \leq 1,\, i, j \in [1,n] \}.\)
For any \(t\in[1,n]\), we define  
\[
\tau_t(M) = \frac{1}{n} \sum_{i=1}^t\lambda_i(M) \quad \quad \text{and}\quad \quad\tau_t(n) = \sup \left\{ \tau_t(M): M \in M_n \right\},
\]  
where \( \lambda_i(M) \)'s are the eigenvalue of \( M \). Furthermore, we define 
\[
\tau_t = \limsup_{n \to \infty} \tau_t(n).
\]  
\begin{lemma}\textnormal{\cite{mohar2009}}\label{sumlargestk}
For $t\geq 1$, we have 
\begin{align*}
    \tau_t&=sup\{\tau_t(M):M\in M_n,n\geq t\}=sup\{\tau_t(M):M\in G_n,n\geq t\}\\
    &=\lim_{n\rightarrow \infty} sup\{\tau_t(M):M\in G_n\},
\end{align*}
where $G_n\subseteq M_n$ is the set of all adjacency matrices of graphs of
order $n$.
\end{lemma}
\begin{theorem}\textnormal{\cite{Chandradas2019}}\label{sumklargestgraph}
    Let $G$ be a graph of order $n$ and $\theta$ be the number of negative eigenvalues of its adjacency matrix. Then for every integer $t\leq n$,
    \[\sum_{i=1}^t\lambda_i\leq \frac{n}{2(\theta+1)}\left(\theta+\sqrt{\theta(\theta t+t-1)}\right).\]
\end{theorem}

We use the classical notation in design theory\cite{stinson2008}. A design is a pair $(X,\mathcal{A)}$, where $X$ is the finite set of elements called points and $\mathcal{A}$ is a collection of non-empty subsets of $X$ called blocks.  Hence, the design naturally corresponds to a hypergraph \( G^* = (V(G^*), E(G^*)) \), where points correspond to vertices, and blocks correspond to hyperedges. 
\begin{definition}\textnormal{\cite{stinson2008}}
    Let $n,\, k$ and $\beta$ be positive integers such that $n>k\geq 2$. A $(n,k,\beta)$-balanced incomplete block design ($(n,k,\beta)$-BIBD) is a design $(X,\mathcal{A})$, where
\begin{itemize}
    \item $\mid X\mid =n$,
    \item each block contains exactly $k$ points, and
    \item every pair of distinct points is contained in exactly $\beta$ blocks.
\end{itemize}
\end{definition}

\begin{theorem}\textnormal{\cite{stinson2008}}\label{BIBDdegree}
    In a $(n,k,\beta)$-BIBD, every point occurs in exactly $\frac{\beta(n-1)}{k-1}$ blocks.
\end{theorem}

A walk of length $s$ on a hypergraph \( G^*\) is a sequence of vertices and edges in $G^*$,  
\(W=v_0e_1v_1\dots\) \( e_sv_s\)
such that \( v_{i-1}, v_i \in e_i \) for $i\in[1,s]$. A subsequence $v_ie_iv_{i+1}\dots e_{j-1}v_j$ of \( W \) is referred to as the \((v_i, v_j)\)-section of \( W \). If all the vertices and edges of $W$ are distinct, then $W$ is called a path. If \( v_0 = v_s \) in a path \( W \), then \( W \) is termed as a cycle. The length of a path or a cycle is its number of edges. Let \( v \in V(G^*) \) and \( e \in E(G^*) \). Then, \( v \)-shrinking on the edge \( e' \) refers to the removal of the vertex \( v \) from the edge \( e' \), resulting in the edge \( e' \setminus \{v\} \).
\begin{table}[H]
    \centering
    \begin{tabular}{l c l}
    \hline
   \\
        $W_s(G^*;u,v)$ & :&The set of $(u,v)$-walks of length $s$ in $G^*$\\
       $M_s(G^*;u,v)$  &:& $|W_s(G^*;u,v)|$\\
       $W_s(G^*;u)$ &:&The set of $(u,u)$-walks of length $s$\\
       $M_s(G^*;u)$ &:&$|W_s(G^*,u)|$ \\
       $(G^*;u)\preceq(G^*;v)$&:& $M_s(G^*;u)\leq M_s(G^*;v)$ for any positive integer $s$\\
       $(G^*;u)\prec(G^*;v)$& :& $(G^*;u)\preceq(G^*;v)$ and $M_{s_0}(G^*;u)< M_{s_0}(G^*;v)$ for \\
       &&some positive integer  $s_0$\\
       \hline
    \end{tabular}
    \caption{Notations}
    \label{tab:my_label}
\end{table}
Let $U(n)$ be the collection of all $k$-uniform unicyclic hypergraphs of order $n$. For $m,k\geq 2,\,\,C_{m,k}=v_1e_1v_2e_2\cdots v_me_mv_1$ be a $k$-uniform cycle of length $m$. For $i\in[1,m]$
\begin{align*}
    e_i&=\{v_i,v_{i+1},u_{i1},u_{i2},\cdots,u_{i(k-2)}\},\, v_{m+1}=v_1,\\
    V_0&=\{v_i,i\in[1,m]\}.
\end{align*}

Let $U(m,n)$ be the class of $k$-uniform unicyclic hypergraph of order $n$ which is constructed from $C_{m,k}$ by attaching $\frac{n}{k-1}-m$ pendant edges to the vertices in $V_0$. Then, $C_m(n_1, n_2,\ldots, n_m)$ refers to a $k$-uniform unicyclic hypergraph obtained from $C_{m,k}$ by attaching the $n_i$ pendant edges to the vertex $v_i$ for each $i\in[1,m]$. In particular, let $X_n \cong C_2(\frac{n}{k-1}-2,0).$  
\begin{example}
    The 3-uniform unicyclic hypergraph $C_3(2,1,0) $ belonging to $U(3,12)$ is illustrated in Figure \ref{fig:1}. 
\end{example}
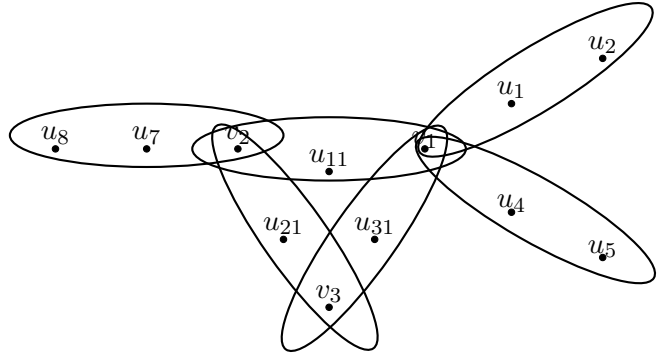
\begin{figure}[H]
\centering

\begin{tikzpicture}[scale=0.6]

 \filldraw[black] (1,0.2) circle (2pt);
    \filldraw[black] (3,-0.3) circle (2pt);
   \filldraw[black] (5.1,0.2) circle (2pt);
   \node(v2) at (1,0.5) {$v_{2}$};
    \node at (3,0) {$u_{11}$};
    \node(v1) at (5.1,0.4) {$v_{1}$};
    \draw[thick] (3,0.2) ellipse (3 and 0.7);
     \filldraw[black] (2,-1.8) circle (2pt);
    \filldraw[black] (3,-3.3) circle (2pt);
   \filldraw[black] (4,-1.8) circle (2pt);
\node(u21) at (2,-1.5) {$u_{21}$};
    \node(v3) at (3,-3) {$v_{3}$};
    \node(u31) at (4,-1.5) {$u_{31}$};
  
    \draw[rotate=35,thick] ($(v2)!0.5!(u21)!0.5!(v3)$) ellipse (0.7 and 3);
     \draw[rotate=55,thick] ($(v1)!0.5!(u31)!0.5!(v3)$) ellipse (3 and 0.7);

    \node(u1) at (7,1.5) {$u_{1}$};
    \filldraw[black] (7,1.2) circle (2pt);
    \node(u2) at (9,2.5) {$u_{2}$};
    \filldraw[black] (9,2.2) circle (2pt);
    \node(u4) at (7,-1) {$u_{4}$};
    \filldraw[black] (7,-1.2) circle (2pt);
    \node(u5) at (9,-2) {$u_{5}$};
    \filldraw[black] (9,-2.2) circle (2pt);
     \draw[rotate=32,thick] ($(v1)!0.5!(u1)!0.5!(u2)$) ellipse (3 and 0.7);
     \draw[rotate=-30,thick] ($(v1)!0.5!(u4)!0.5!(u5)$) ellipse (3 and 0.7);

       \node(u7) at (-1,0.5) {$u_{7}$};
       \node(u8) at (-3,0.5) {$u_{8}$};
        \filldraw[black] (-1,0.2) circle (2pt);
         \filldraw[black] (-3,0.2) circle (2pt);
          \draw[thick] (-1,0.5) ellipse (3 and 0.7);
\end{tikzpicture}
\caption{Hypergraph $C_3(2,1,0)$ }
\label{fig:1}
\end{figure}

\begin{lemma}\textnormal{\cite{feng2023}}\label{geq3}
Let $G_1^*\in U(m,n)$ be a $k$-uniform hypergraph with a unique cycle $C_{m,k}$ and $n_i$ be the number of pendant edges attached to each vertex $v_i\in V_0$. Let $G_2^*$ be a hypergraph obtained from $G_1^*$ by $v_1$-shrinking on $e_m$, and $e_0=e_m\setminus \{v_1\}$. If $m\geq 4$ and $n_3=max_{i\in [1,m]}n_i$, then $EE(G_1^*)<EE(G_2^*-e_0+(e_0\cup \{v_3\})).$\\
Note that, $G_2^*-e_0+(e_0\cup \{v_3\})$ is a $k$-uniform unicyclic hypergraph of order $n$ and the length of its unique cycle is $m-2$.
\end{lemma}
\begin{lemma}\textnormal{\cite{feng2023}}\label{C3ni}
    If $n_1\geq n_2 \geq n_3 \geq 1,$ then \[EE(C_3(n_1,n_2,n_3))<EE(C_3(n_1+n_2,n_3,0))<EE(C_3(n_1+n_2+n_3,0,0)).\]
\end{lemma}
\begin{lemma}\textnormal{\cite{Wand2020}}\label{EEuv}
    For a hypergraph $G^*$  with $u,v\in V(G^*)$ and 
    $E_0=\{e_1,e_2,\cdots,e_s\}$, where $e_i\in \mathcal{P}(V(G^*)),\, e_i\cap \{u,v\}=\phi$ for $i\in [1,s]$. Let 
    \begin{align*}
        E_u&=\{e_i'=e_i\cup\{u\}:e_i\in E_0,i\in [1,s]\},\\
         E_v&=\{e_i'=e_i\cup\{v\}:e_i\in E_0,i\in [1,s]\},
    \end{align*}
    $G^*_u=G^*-E_0+E_u$ and $G^*_v=G^*-E_0+E_v$. If for any $w\in \bigcup_{i\in[1,s]}e_i,\,(G^*,u)\prec (G^*,v)$ and $(G^*:w,u)\preceq (G^*:w,v)$, then $EE(G^*_u)<EE(G^*_v)$.
\end{lemma}
\begin{lemma}\textnormal{\cite{Wand2020}}\label{EEcolescence}
    Let $G^*_1$ and $G^*_2$ be two hypergraphs with $u,v \in V(G^*_1)$ and $w\in V(G^*_2)$. If $(G^*_1;u)\prec (G^*_1;v)$, then $EE(G^*_1(u)\cdot (w)G^*_2)< EE(G^*_1(v)\cdot (w)G^*_2)$, where $G^*_1(u)\cdot (w)G^*_2$ is the hypergraph obtained from $G^*_1$ and $G^*_2$ by identifying $u$ of $G^*_1$ with $w$ of $G^*_2$.
\end{lemma}
\begin{lemma}\textnormal{\cite{feng2023}}\label{numberofwalk}
    For a $k$-uniform unicyclic hypergraph $G^*$ with a unique cycle $C_{m,k},$ let $u \in V_0$ and $u_1$ be a pendant neighbour of $u$ in $G^*$, then $(G^*;u_1)\prec (G^*;u)$.
\end{lemma}
\begin{theorem}\textnormal{\cite{feng2023}}\label{unicyclicmax}
  Let $G^*$ be a $k$-uniform unicyclic hypergraph of order $n$, $EE(G^*)\leq EE(C_2(\frac{n}{k-1}-2,0))$ and the equality holds if and only if $G^* \cong C_2(\frac{n}{k-1}-2,0)$.
\end{theorem}

\begin{theorem}\textnormal{\cite{kurian2024}}\label{sumlamda2}
     Let $G^*$ be a $k$- uniform hypergraph with $n$ vertices and $m$ edges, and $\lambda_i,\,i\in[1,n]$ be the adjacency eigenvalues of $G^*$. Then
    \begin{equation*}
        k(k-1)m\leq\sum_{i=1}^n\lambda_i^2\leq (k-1)m\Bigl(m(k-2)+2\Bigr).
    \end{equation*}
\end{theorem}

\section{Eigenvalues of $k$-uniform hypergraph}
In this section, we study some of the interesting properties of the eigenvalues of a $k$-uniform hypergraph. 
\begin{theorem}
    Let $B$ be a symmetric matrix of order n whose largest entry is $b$ and the smallest entry is $a$. Then for every integer $t,\,2\leq t\leq n$,
    \[\tau_t \leq \dfrac{1}{2\theta+1}\left(\theta+\sqrt{\theta(t\theta +t-1)}\right)(b-a)+max\{0,a\},\]
where $\theta$ is the number of negative eigenvalues of $B$.
\end{theorem}
\begin{proof}
From Lemma \ref{sumlargestk}, to determine $\tau_t$ of symmetric matrices it is enough to consider the adjacency matrices of graphs. Therefore, from Theorem \ref{sumklargestgraph}
\[ \tau_t(M)\leq \dfrac{1}{2\theta+1}\left(\theta+\sqrt{\theta(t\theta +t-1)}\right),\,M\in M_n.\]
Let $C=B-aJ_n$. Then, $\dfrac{1}{b-a}C\in M_n$. Hence,
\begin{align*}
    \frac{1}{(b-a)n}\sum_{i=1}^t\lambda_i(C)&\leq\dfrac{1}{2\theta+1}\left(\theta+\sqrt{\theta(t\theta +t-1)}\right),\\
    \sum_{i=1}^t\lambda_i(C)&\leq\dfrac{n(b-a)}{2\theta+1}\left(\theta+\sqrt{\theta(t\theta +t-1)}\right).
\end{align*}
Since \[\sum_{i=1}^t\lambda_i(B)\leq \sum_{i=1}^t\lambda_i(C)+\sum_{i=1}^t\lambda_i(aJ_n)\] and \[\sum_{i=1}^t\lambda_i(aJ_n)=\begin{cases}
    an &if \quad  a\geq 0,\\
    0  &if \quad a<0,
\end{cases}\]  we get
\[\sum_{i=1}^t\lambda_i(B)\leq\dfrac{n(b-a)}{2\theta+1}\left(\theta+\sqrt{\theta(t\theta +t-1)}\right)+max\{0,an\}.\]
Consequently, this establishes the result.

\end{proof}
\begin{corollary}
        Let $G^*$ be a $k$-uniform hypergraph of order $n$. Then for every integer $t,\,t\in[2,n]$,
    \[\sum_{i=1}^t\lambda_i \leq \dfrac{n\binom{n-2}{k-2}}{2\theta+1}\left(\theta+\sqrt{\theta(t\theta +t-1)}\right),\]
where $\theta$ is the number of negative eigenvalues of $A(G^*)$.
\end{corollary}

\begin{theorem}\label{distincteig}
    Let $G^*$ be a $k$-uniform hypergraph on $n$ vertices. Then $G^*=(V(G^*),E(G^*))$ has 2 distinct eigenvalues if and only if $G^*$ is a $(n,k,\beta)$- BIBD, where every pair of vertices is contained in exactly $\beta$ edges.
\end{theorem}
\begin{proof}
    Suppose that $G^*$ has only two distinct eigenvalues. By the Perron-Frobenius theorem, $\lambda_1$ and $-\dfrac{\lambda_1}{n-1}$ are the eigenvalues of $G^*$ with multiplicity 1 and $(n-1)$, respectively. Therefore,
    \[\text{Rank}\left(A(G^*)+\frac{\lambda_1}{n-1}I_n\right)=1.\]
     Since every rank-one matrix can be represented as the outer product of two column vectors, there exist two vectors  \(\mathbf{v} = \begin{bmatrix} v_1 & v_2 & v_3 & \cdots & v_n \end{bmatrix}^T\) 
and  
\(
\mathbf{w} = \begin{bmatrix} w_1 & w_2 & w_3 & \cdots & w_n \end{bmatrix}^T
\)such that  
\[
A(G^*)+\frac{\lambda_1}{n-1}I_n = \mathbf{v w}^T =  
\begin{bmatrix}  
    v_1w_1 & v_1w_2 & v_1w_3 & \cdots & v_1w_n \\  
    v_2w_1 & v_2w_2 & v_2w_3 & \cdots & v_2w_n \\  
    \vdots & \vdots & \vdots & \ddots & \vdots \\  
    v_nw_1 & v_nw_2 & v_nw_3 & \cdots & v_nw_n  
\end{bmatrix}.
\]
Also, since $A(G^*)$ is symmetric and its diagonal entries are all zero, the following conditions must be satisfied:
\begin{itemize}
    \item $v_iw_i=\dfrac{\lambda_1}{n-1},i\in[1,n] $,
    \item $v_iw_j=v_jw_i,\,i\neq j,\,i,j\in[1,n].$
\end{itemize}
On solving these equations, we get $A(G^*)=\frac{\lambda_1}{n-1}(J_n-I_n)=\beta(J_n-I_n)$. Thus, every two vertices is contained in exactly $\beta$ edges. Hence, $G^*$ is an $(n,k,\beta)$-BIBD. The converse holds trivially. 
\end{proof}
The following corollary is an immediate consequence of Theorems \ref{BIBDdegree} and \ref{distincteig}.
\begin{corollary}
    Let $G^*$ be a $k$-uniform hypergraph of order $n$ with exactly two distinct eigenvalues. Then $G^*$ is regular with
    \begin{equation*}
        d(v_i)=\frac{\beta(n-1)}{k-1}, i\in[1,n],
    \end{equation*}
    where every pair of vertices is contained in exactly $\beta$ edges.
\end{corollary}

\section{Estrada index of a $k$-uniform hypergraph}
In \cite{Rodriguez2021}, the authors establish a lower bound for the Estrada index of graphs. Applying similar arguments, we establish that for a hypergraph $G^*$ of order $n$, a lower bound for its Estrada index is given by 
\[EE(G^*)\geq e^{\lambda_1}+(n-1)-\lambda_1.\]
 Now, let $f(x)=e^x+n-1-x$ be a function on $x\geq 0$. Since $f(x)$ increases, any lower bound for $\lambda_1$ provides another lower bound for the Estrada index of $G^*$. The following theorem establishes a lower bound for the Estrada index in terms of order, number of edges, and uniformity of the hypergraph.
\begin{theorem}
    Let $G^*$ be a $k$-uniform hypergraph with $n$ vertices and $m$ edges. Then \[EE(G^*)\geq \sqrt{n^2+\dfrac{4k(k-1)m}{2}}.\]
     The equality holds if and only if $G^*=\overline{K}_n^k$.
\end{theorem}
\begin{proof}
    Note that
    \begin{equation}\label{eqn1}
        (EE(G^*))^2=\sum_{i=1}^n e^{2\lambda_i}+2\sum_{i<j}e^{\lambda_i}e^{\lambda_j}.
    \end{equation}
    Using AM-GM inequality, we get
    \begin{align}\label{eqn2}
        \dfrac{2\sum_{i<j}e^{\lambda_i}e^{\lambda_j}}{n(n-1)}&\geq \left(\prod_{i<j}e^{\lambda_i}e^{\lambda_j}\right)^{\frac{2}{n(n-1)}},\nonumber\\
        2\sum_{i<j}e^{\lambda_i}e^{\lambda_j}&\geq n(n-1)\left(\prod_{i=1}^ne^{\lambda_i}\right)^{\frac{2}{n}}=n(n-1).
    \end{align}
Using power series expansion and spectral moments, we have
\begin{align*}
    \sum_{i=1}^ne^{2\lambda_i}&=\sum_{i=1}^n\sum_{s\geq 0} \dfrac{(2\lambda_i)^s}{s!}=n+2\sum_{i=1}^n\lambda_i^2+\sum_{i=1}^n\sum_{s\geq3}\dfrac{(2\lambda_i)^s}{s!}.
\end{align*}
Applying Theorem \ref{sumlamda2},
\[\sum_{i=1}^ne^{2\lambda_i}\geq n+2k(k-1)m+\sum_{i=1}^n\sum_{s\geq3}\dfrac{(2\lambda_i)^s}{s!}.\]
Aiming at a lower bound we replace $\sum_{s\geq3}\dfrac{(2\lambda_i)^s}{s!}$ by $8\sum_{s\geq3}\dfrac{(\lambda_i)^s}{s!}$. Therefore, 
\begin{align}\label{eqn3}
    \sum_{i=1}^ne^{2\lambda_i}&\geq n+2k(k-1)m+\gamma \sum_{i=1}^n\sum_{s\geq 3}\dfrac{(\lambda_i)^s}{s!}\nonumber\\
    &\geq n+2k(k-1)m-\gamma n-\gamma (k-1)m\frac{(m(k-2)+2)}{2}+\gamma \sum_{i=1}^n\sum_{s\geq 0} \dfrac{(\lambda_i)^s}{s!}\nonumber\\
    &=(1-\gamma)n+(2k-\dfrac{\gamma(m(k-2)+2)}{2})(k-1)m+\gamma EE(G^*), \,\qquad\gamma\in[0,8].
\end{align}
From (\ref{eqn1}),(\ref{eqn2}) and (\ref{eqn3}), we have
\[(EE(G^*))^2-\gamma EE(G^*)-n(n-\gamma)+(4k-\gamma(m(k-2)+2))\frac{(k-1)m}{2}\geq 0.\]
On solving, we get
\[EE(G^*)\geq \dfrac{\gamma}{2}+\sqrt{\left(\dfrac{\gamma}{2}-n\right)^2+(4k-\gamma(m(k-2)+2))\dfrac{(k-1)m}{2}}.\]
Now, consider the function
\[f(x)= \dfrac{x}{2}+\sqrt{\left(\dfrac{x}{2}-n\right)^2+(4k-x(m(k-2)+2))\dfrac{(k-1)m}{2}}\]
which is monotonically decreasing in [0,8]. Therefore, the best lower bound is at $\gamma=0$. Hence
\[EE(G^*)\geq \sqrt{n^2+\dfrac{4k(k-1)m}{2}}.\]
Clearly, equality holds if and only if the adjacency matrix of $G^*$ has no non-zero eigenvalue.  That is $G^*$ is an edgeless hypergraph $\overline{K}_n^k$.
\end{proof}
\begin{theorem}\label{thm4.2}
 Let $G^*$ be a $k$-uniform hypergraph with $n$ vertices and $m$ edges. Then \[EE(G^*)\leq n-1+e^{\sqrt{(k-1)m(m(k-2)+2)}}. \]   
  The equality holds if and only if $G^*=\overline{K}_n^k$.
\end{theorem}
\begin{proof}
    From the definition of Estrada index and Theorem \ref{sumlamda2}, we have
    \begin{align*}
        EE(G^*)&=n+\sum_{i=1}^n\sum_{s\geq 1}\dfrac{(\lambda_i)^s}{s!}\\
        &\leq n+\sum_{i=1}^n\sum_{s\geq 1}\dfrac{|\lambda_i|^s}{s!} \leq n+\sum_{s\geq 1}\dfrac{1}{s!}\left(\sum_{i=1}^n(\lambda_i^2)\right)^\frac{s}{2}\\
        &\leq n+\sum_{s\geq 1}\dfrac{(\sqrt{
        (k-1)m(m(k-2)+2)})^s}{s!}=n-1+e^{\sqrt{(k-1)m(m(k-2)+2)}}.
    \end{align*}
    The equality holds if and only if the adjacency matrix of $G^*$ has no non-zero eigenvalue.  That is, $G^*=\overline{K}_n^k$.
\end{proof}
The energy of a hypergraph $G^*$is defined as $E(G^*)=\sum_{i=1}^n \mid \lambda_i \mid.$
In the following theorem, we obtain a lower bound of the Estrada index of the $k$-uniform hypergraph in terms of its energy.  
\begin{theorem}
    Let $G^*$ be a $k$-uniform hypergraph with $n$ vertices and $m$ edges, and let $E(G^*)$ be its energy. Then
    \[EE(G^*)\leq n+E(G^*)-1-\sqrt{(k-1)m\left(m(k-2)+2\right)}+e^{\sqrt{(k-1)m\left(m(k-2)+2\right)}}.\]
    The equality holds if and only if $G^*=\overline{K}_n^k$.
    
\end{theorem}
\begin{proof}
    Using similar arguments from the proof of the Theorem \ref{thm4.2}, we get
    \begin{align*}
        EE(G^*)&\leq n+E(G^*)+\sum_{s\geq2}\frac{1}{s!}\left(\sqrt{(k-1)m\left(m(k-2)+2\right)}\right)^s\\
        &=n+E(G^*)-1-\sqrt{(k-1)m\left(m(k-2)+2\right)}+e^{\sqrt{(k-1)m\left(m(k-2)+2\right)}}.
    \end{align*}
\end{proof}
\begin{remark}
    From the definition of the Estrada index of hypergraph $G^*$ of order $n$, we have
    \begin{align*}
        EE(G^*)\leq n+\sum_{i=1}^n\sum_{s\geq1}\frac{\mid\lambda_i\mid^s}{s!}=n+\sum_{s\geq1}\frac{(E(G^*))^s}{s!}.
    \end{align*}
    Therefore, \[EE(G^*)\leq n-1+e^{E(G^*)}.\]
    and the equality holds if and only if $G^*=\overline{K}_n^k.$
\end{remark}

Next, we provide the Nordhaus–Gaddum-type bounds for the Estrada index of a $k$-uniform hypergraph.
\begin{theorem}
    Let $G^*$ be a $k$-uniform hypergraph of order $n$ and $\overline{G^*}$  be its complement. Then 
    \[EE(G^*)+EE(\overline{G^*})\geq 2 e^{\frac{n-1}{2}}+2(n-1)e^{-\frac{1}{2}}.\]
    Equality holds if and only if $G^*$ is a $k$-uniform hypergraph satisfying $(n,k,\beta)$-BIBD.
\end{theorem}
\begin{proof}
    Applying the AM-GM inequality, we get
    \begin{align*}
        EE(G^*)&
        &\geq e^{\lambda_1}+(n-1)e^{-\frac{\lambda_1}{(n-1)}}.
    \end{align*}
    Let $f(x)=e^x+(n-1)e^{-\frac{x}{(n-1)}}$.
    This function is increasing for all $x\geq 0$. Let $\lambda_i$ be the eigenvalues of $G^*$, and $\overline{\lambda}_i$ be the eigenvalues of $\overline{G^*}$ arranged in the decreasing order $\overline{\lambda}_1\geq \overline{\lambda}_2\geq \overline{\lambda}_3\ldots,\overline{\lambda}_n$. Then
    \begin{align*}
        EE(G^*)+EE(\overline{G^*})&\geq e^{\lambda_1}+e^{\overline{\lambda}_1}+(n-1)\left( e^{-\frac{\lambda_1}{(n-1)}}+e^{-\frac{\overline{\lambda}_1}{(n-1)}} \right)\\
        &\geq 2e^{\frac{\lambda_1+\overline{\lambda}_1}{2}}+2(n-1)e^{-\frac{\lambda_1+\overline{\lambda}_1}{n-1}}.
    \end{align*}
    Since $\lambda_1+\overline{\lambda}_1\geq n-1$, we have
    \[  EE(G^*)+EE(\overline{G^*})\geq 2 e^{\frac{n-1}{2}}+2(n-1)e^{-\frac{1}{2}}.\]
   Equality holds if and only if $G^*$ is a $k$-uniform hypergraph with two distinct eigenvalues, say $\lambda_1$ and $-\frac{\lambda_1}{n-1}$. Hence, by Theorem \ref{distincteig} $ G^*$ is a $k$-uniform hypergraph that satisfies $(n,k,\beta)$-BIBD.
\end{proof}

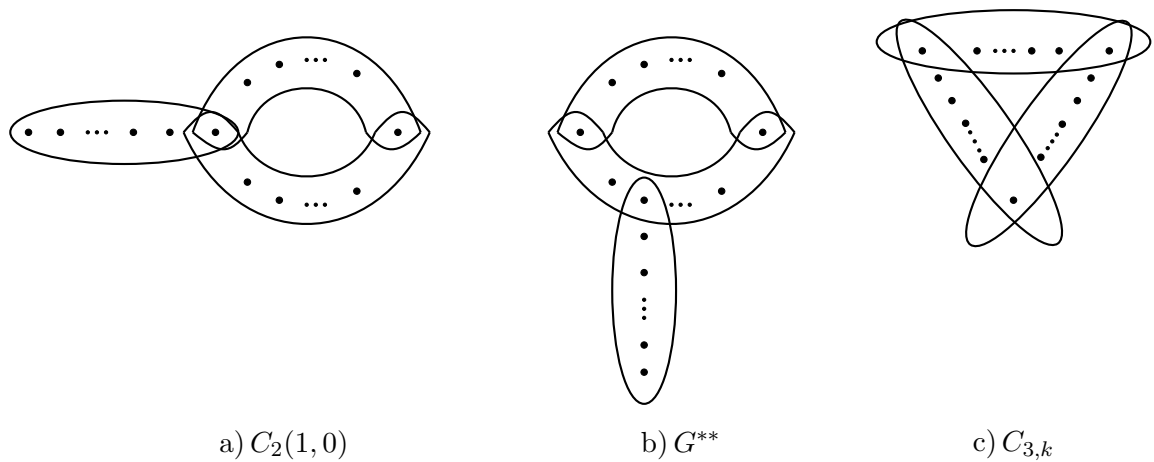
\begin{figure}[H]
    \centering
\begin{tikzpicture}[scale=0.6]

 \filldraw[black] (1,0) circle (2pt);
    \filldraw[black] (2.2,0) circle (2pt);
    \filldraw[black] (3.4,0) circle (2pt);
    \filldraw[black] (4,0) circle (2pt);
    \filldraw[black] (2.6,0) circle (1pt);
    \filldraw[black] (2.8,0) circle (1pt);
    \filldraw[black] (3.0,0) circle (1pt);
   \filldraw[black] (5.1,0.0) circle (2pt);
   \node(v2) at (1,0.3) {};
    \node at (3,0) {};
    \node(v1) at (5.1,0.2) {};
    \draw[thick] (3,0.2) ellipse (3 and 0.7);
    \filldraw[black] (1.35,-0.6) circle (2pt);
    \filldraw[black] (1.65,-1.1) circle (2pt);
    \filldraw[black] (1.95,-1.6) circle (2pt);
    \filldraw[black] (2.05,-1.8) circle (1pt);
    \filldraw[black] (2.15,-2) circle (1pt);
    \filldraw[black] (2.25,-2.2) circle (1pt);
    \filldraw[black] (2.35,-2.4) circle (2pt);
     \filldraw[black] (3,-3.3) circle (2pt);
    \filldraw[black] (3.6,-2.35) circle (2pt);
    \filldraw[black] (3.75,-2.15) circle (1pt);
    \filldraw[black] (3.85,-2.0) circle (1pt);
    \filldraw[black] (3.95,-1.85) circle (1pt);
    \filldraw[black] (4.07,-1.62) circle (2pt);
    \filldraw[black] (4.4,-1.1) circle (2pt);
    \filldraw[black] (4.7,-0.6) circle (2pt);
\node(u21) at (2,-1.5) {};
    \node(v3) at (3,-3) {};
    \node(u31) at (4,-1.5) {};
  
    \draw[rotate=35,thick] ($(v2)!0.5!(u21)!0.5!(v3)$) ellipse (0.7 and 3);
     \draw[rotate=55,thick] ($(v1)!0.5!(u31)!0.5!(v3)$) ellipse (3 and 0.7);
\node at (3,-8.7) {$\mathrm{c)}\,C_{3,k}$};
\draw[thick] (-7,-1.8) .. controls (-6,1) and (-3,1) .. (-2,-1.8);
\draw[thick] (-5.8,-1.8) .. controls (-5.5,-0.5) and (-3.5,-0.5) .. (-3.2,-1.8);
\draw[thick] (-5.8,-1.8) .. controls (-6.2,-2.3) and (-6.4,-2.3)  .. (-7,-1.8);
\draw[thick] (-2,-1.8) .. controls (-2.5,-2.3) and (-2.8,-2.3) .. (-3.2,-1.8);

\draw[thick] (-7.2,-1.8) .. controls (-6,-4.5) and (-3,-4.5) .. (-1.8,-1.8);
\draw[thick] (-6,-1.8) .. controls (-5.5,-3.1) and (-3.5,-3.1) .. (-3,-1.8);
\draw[thick] (-6,-1.8) .. controls (-6.5,-1.2) and (-6.7,-1.2)  .. (-7.2,-1.8);
\draw[thick] (-1.8,-1.8) .. controls (-2.3,-1.2) and (-2.5,-1.2) .. (-3,-1.8);
\filldraw[black] (-6.5,-1.8) circle (2pt);
\filldraw[black] (-2.5,-1.8) circle (2pt);
\filldraw[black] (-5.8,-2.9) circle (2pt);
\filldraw[black] (-5.1,-3.3) circle (2pt);
\filldraw[black] (-4.5,-3.4) circle (1pt);
\filldraw[black] (-4.3,-3.4) circle (1pt);
\filldraw[black] (-4.1,-3.4) circle (1pt);
\filldraw[black] (-3.4,-3.1) circle (2pt);

\filldraw[black] (-5.8,-.7) circle (2pt);
\filldraw[black] (-5.1,-.3) circle (2pt);
\filldraw[black] (-4.5,-0.2) circle (1pt);
\filldraw[black] (-4.3,-0.2) circle (1pt);
\filldraw[black] (-4.1,-0.2) circle (1pt);
\filldraw[black] (-3.4,-0.5) circle (2pt);

 \draw[thick] (-5.1,-5.3) ellipse (0.7 and 2.5);
\filldraw[black] (-5.1,-4.1) circle (2pt);
\filldraw[black] (-5.1,-4.9) circle (2pt);
\filldraw[black] (-5.1,-5.5) circle (1pt);
\filldraw[black] (-5.1,-5.7) circle (1pt);
\filldraw[black] (-5.1,-5.9) circle (1pt);
\filldraw[black] (-5.1,-6.5) circle (2pt);
\filldraw[black] (-5.1,-7.1) circle (2pt);
\node at (-4.3,-8.7) {$\mathrm{b)}\,G^{**}$};
\draw[thick] (-15,-1.8) .. controls (-14,1) and (-11,1) .. (-10,-1.8);
\draw[thick] (-13.8,-1.8) .. controls (-13.5,-0.5) and (-11.5,-0.5) .. (-11.2,-1.8);
\draw[thick] (-13.8,-1.8) .. controls (-14.2,-2.3) and (-14.4,-2.3)  .. (-15,-1.8);
\draw[thick] (-10,-1.8) .. controls (-10.5,-2.3) and (-10.8,-2.3) .. (-11.2,-1.8);

\draw[thick] (-15.2,-1.8) .. controls (-14,-4.5) and (-11,-4.5) .. (-9.8,-1.8);
\draw[thick] (-14,-1.8) .. controls (-13.5,-3.1) and (-11.5,-3.1) .. (-11,-1.8);
\draw[thick] (-14,-1.8) .. controls (-14.5,-1.2) and (-14.7,-1.2)  .. (-15.2,-1.8);
\draw[thick] (-9.8,-1.8) .. controls (-10.3,-1.2) and (-10.5,-1.2) .. (-11,-1.8);
\filldraw[black] (-14.5,-1.8) circle (2pt);
\filldraw[black] (-10.5,-1.8) circle (2pt);
\filldraw[black] (-13.8,-2.9) circle (2pt);
\filldraw[black] (-13.1,-3.3) circle (2pt);
\filldraw[black] (-12.5,-3.4) circle (1pt);
\filldraw[black] (-12.3,-3.4) circle (1pt);
\filldraw[black] (-12.1,-3.4) circle (1pt);
\filldraw[black] (-11.4,-3.1) circle (2pt);

\filldraw[black] (-13.8,-.7) circle (2pt);
\filldraw[black] (-13.1,-.3) circle (2pt);
\filldraw[black] (-12.5,-0.2) circle (1pt);
\filldraw[black] (-12.3,-0.2) circle (1pt);
\filldraw[black] (-12.1,-0.2) circle (1pt);
\filldraw[black] (-11.4,-0.5) circle (2pt);

 \draw[thick] (-16.5,-1.8) ellipse (2.5 and 0.7);
\filldraw[black] (-15.5,-1.8) circle (2pt);
\filldraw[black] (-16.3,-1.8) circle (2pt);
\filldraw[black] (-16.9,-1.8) circle (1pt);
\filldraw[black] (-17.1,-1.8) circle (1pt);
\filldraw[black] (-17.3,-1.8) circle (1pt);
\filldraw[black] (-17.9,-1.8) circle (2pt);
\filldraw[black] (-18.6,-1.8) circle (2pt);
\node at (-13,-8.7) {$\mathrm{a)}\,C_2(1,0)$};
\end{tikzpicture}
\caption{$k$-uniform unicyclic hypergraph when $n=3(k-1)$}
\label{fig.3}
\end{figure}

From (\ref{eqn1a}), it is evident that for any hypergraph $G^*$ of order $n$, adding an edge $e'$ contributes to an increase in the number of closed walks. Consequently, this leads to an increase in the Estrada index of the hypergraph. Therefore,
\[EE(G^*)<EE(G^*+e').\]
Also, it is evident that for a $k$- uniform hypergraph of order $n$,
\[EE(\overline{K}_n^k)<EE(G^*)<EE(K_n^k).\]
From Theorem \ref{unicyclicmax}, we have $X_n$ has the maximum Estrada index among the $k$-uniform unicyclic hypergraph. Now, we obtain the $k$-uniform unicyclic hypergraph with the second largest Estrada index. When \( n = 3(k - 1) \), there exist only three types of $k$-uniform unicyclic hypergraph, as illustrated in Figure \ref{fig.3}. Among these, \( C_2(1,0) \) possesses the maximum Estrada index and as discussed later in Remark~\ref{remark4.11} it is clear that the $k$-uniform hypergraph \( G^{**} \)($k$-uniform hypergraph with a cycle of length 2 and an additional edge incident to one of the pendant vertices cycle) has the second largest Estrada index. When \( n < 3(k - 1) \), the only possible structure is the cycle \( C_{2,k} \). Therefore, the following result focuses on the case where \( n > 3(k - 1) \).

\begin{theorem}\label{secmaxEEG}
    Let $G^{*}_0$ in $U(n) \setminus \{X_n\},\, n>3(k-1)$ be the $k$-uniform hypergraph with maximum Estrada index. Then $G^{*}_0\in U(m,n)\setminus\{X_n\}.$
\end{theorem}
\begin{proof}
Assume that $G^*_0$ in $U(n) \setminus \{X_n\},\, n>3(k-1)$ is the $k$-uniform hypergraph with the maximum Estrada index such that $G^*_0\notin U(m,n)\setminus\{X_n\}.$ \\
    \textbf{Case I.} When $m\geq 3$. Since $G^*_0\notin U(m,n)\setminus\{X_n\}$, there exists at least one vertex $v_t\in V_0$ with a non-pendant neighbour $u_t$ outside $C_{m,k}$ or in $V(C_{m,k})\setminus V_0$. Let 
    \begin{align*}
        E_0&=\{\Tilde{e}_i=e_i \setminus \{u_t\}:e_i \text{ is incident to } u_t \text{ and }  v_t\notin e_i\},\\
        &=\{\Tilde{e}_1,\Tilde{e}_2,\cdots,\Tilde{e}_{d(u_t)-1}\},\\
        E_{v_t}&=\{ \Tilde{e}_i\cup \{v_t\}:\Tilde{e}_i\in E_0,i\in[1,d(u_t)-1]\},\\
        E_{u_t}&=\{ \Tilde{e}_i\cup \{u_t\}:\Tilde{e}_i\in E_0,i\in[1,d(u_t)-1]\}.
    \end{align*}
In $G^*_0-E_{u_t}$, we have the following
\begin{itemize}
    \item $d(u_t)=1$ and $d(v_t)\geq 2$, then by Lemma \ref{numberofwalk} $(G^*_0-E_{u_t};u_t)\prec (G^*_0-E_{u_t};v_t)$.
    \item For any $w\in E_0,\,M_k(G^*_0-E_{u_t};w,u_t)=M_k(G^*_0-E_{u_t};w,v_t)=0$.
\end{itemize}
By Lemma \ref{EEuv}, we get $EE(G^*_0)<EE(G^*_0-E_{u_t}+E_{v_t})$. This is a contradiction of our assumption. Hence $G^*_0\in U(m,n)\setminus\{X_n\}.$ \\
\textbf{Case II.}When $m=2$. Since $G^*_0\notin U(m,n)\setminus\{X_n\}$ and $n > 3(k-1)$,  $G^*_0$ must satisfy one of the following conditions:
\begin{itemize}
    \item[i.] Both vertices in $V_0$ are incident to an edge outside  $C_{2,k}$, or one vertex in  $V_0$ and one in $V(C_{2,k})\setminus V_0$  are incident to an edge outside  $C_{2,k}$.
    \item[ii.] Exactly one vertex in $V_0$ is adjacent to a non-pendant vertex either outside $C_{2,k}$ or in $V(C_{2,k})\setminus V_0$, and at least one such non-pendant vertex is adjacent to another non-pendant vertex distinct from vertices of $V_0$ 
    \item[iii.] Exactly one vertex in $V_0$ is adjacent to a non-pendant vertex outside $C_{2,k}$ and the neighbours of this non-pendant vertex are all pendant vertices except for the vertex in $V_0$.
\end{itemize}
  The conditions i. and ii. contradict our assumption by following the same arguments as in Case I. Now we only need to prove for $G^*_0$ satisfying condition iii. Without loss of generality, $v_2$ is adjacent to a non-pendant vertex $u_t$ outside $C_{2,k}$ and $v_1$ does not have any neighbour outside $C_{2,k}$.
 Let $\bar{E}_0=\{e: e\in E(G^*_0)$ be an incident edge to $u_t$, $v_2\notin e\}.$ The edges of $\bar{E}_0$ form a k-uniform hyperstar $S^k_*$. Take $G^*_1=G^*_0-V(S^k_*)\setminus\{u_t\}.$ Clearly, we have $(G^*_1,u_t)\prec(G^*_1,v_1)$. By Lemma \ref{EEcolescence}, We get 
 $EE(G^*_0)<EE(G^*_1(v_1)\cdot (u_t)S^k_*)$, a contradiction. Hence $G^*_0\in U(m,n)\setminus\{X_n\}$.
\end{proof}
\begin{lemma}\label{lemma 4.2}
    If $n_1\geq n_2\geq 1$, then $EE(C_2(n_1,n_2))<EE(C_2(n_1+1,n_2-1))
    $.
\end{lemma}
\begin{proof}
    Let $e_{21}$ be one of the edges attached to the vertex $v_2\in V_0$ and $G^*=C_2(n_1,n_2)-e_{21}$. Then, clearly $(G^*;v_2)\prec(G^*;v_1)$ and for $w\in e_{21}\setminus \{v_2\}$, $M_s(G^*;w,v_2)=M_s(G^*;w,v_1)=0$. Then, the result follows from Lemma \ref{EEuv}.
\end{proof}
\begin{lemma}\label{lemma4.3}
    $EE(C_3(\dfrac{n}{k-1}-3,0,0))<EE(C_2(\dfrac{n}{k-1}-3,1))$.
\end{lemma}
\begin{proof}
    Let $G'^*$ be the hypergraph obtained from $C_3(\dfrac{n}{k-1}-3,0,0)$ by $v_3$ shrinking on $e_3$ and $G$ be the graph obtained from $G'^*$ by replacing each edge of $G'^*$ by a complete graph $K_{|e|}$ (Figure \ref{fig:2} illustrates the hypergraph $G'^*$ and the corresponding graph $G$ obtained from a 3-uniform hypergraph $C_3(2,0,0)$ on 10 vertices).\\
    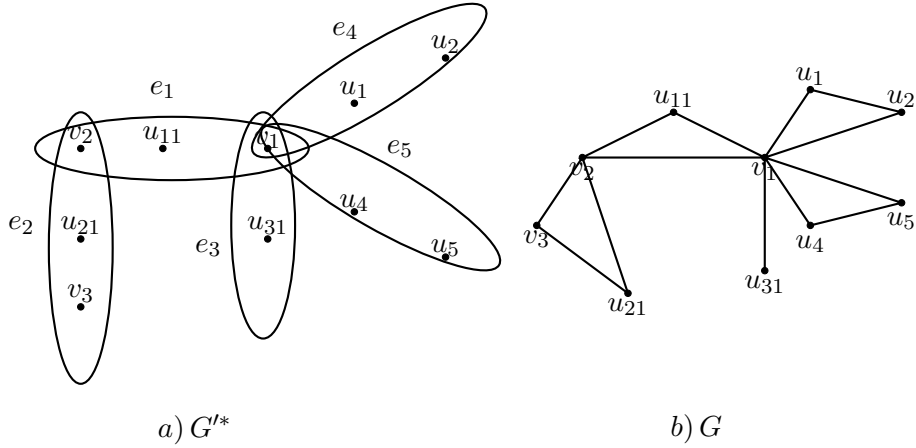
\begin{figure}[H]
\centering
\begin{tikzpicture}[scale=0.6]
 \filldraw[black] (1,0.2) circle (2pt);
   \filldraw[black] (5.1,0.2) circle (2pt);

   \filldraw[black] (2.8,0.2) circle (2pt);
   \node at (1,0.5) {$v_{2}$};
     \node at (2.8,0.5) {$u_{11}$};
     \node at (2.8,1.5) {$e_{1}$};
    \node(v1) at (4.3,0.3) {};
     \node at (5.1,0.4) {$v_{1}$};  
    \draw[thick] (3,0.2) ellipse (3 and 0.7);
     \filldraw[black] (1,-1.8) circle (2pt);
    \filldraw[black] (1,-3.3) circle (2pt);
   \filldraw[black] (5.1,-1.8) circle (2pt);
  
\node at (1,-1.5) {$u_{21}$};
\node at (-0.3,-1.5) {$e_{2}$};
    \node at (1,-3) {$v_{3}$};
    \node at (5.1,-1.5) {$u_{31}$};
    \node at (3.8,-2) {$e_3$};
     \draw[thick] (1,-2) ellipse (0.7 and 3);
      \draw[thick] (5,-1.5) ellipse (0.7 and 2.5);
      
    \node(u1) at (7,1.5) {$u_{1}$};
     \node at (6.8,2.8) {$e_4$};
    \filldraw[black] (7,1.2) circle (2pt);
    \node(u2) at (9,2.5) {$u_{2}$};
    \filldraw[black] (9,2.2) circle (2pt);
    \node(u4) at (7,-1) {$u_{4}$};
    \node(u4) at (8,0.2) {$e_5$};
    \filldraw[black] (7,-1.2) circle (2pt);
    \node(u5) at (9,-2) {$u_{5}$};
    \filldraw[black] (9,-2.2) circle (2pt);
     \draw[rotate=32,thick] ($(v1)!0.5!(u1)!0.5!(u2)$) ellipse (3 and 0.7);
     \draw[rotate=-30,thick] ($(v1)!0.5!(u4)!0.5!(u5)$) ellipse (3 and 0.7);
 \node(u11) at (14,1.3) {$u_{11}$};
 \node(v2) at (12,-0.3) {$v_{2}$};
 \node(v1) at (16,-0.3) {$v_{1}$};

\draw[thick] (14,1) -- (12,0)--(16,0)--(14,1);
\filldraw[black] (14,1) circle (2pt);
\filldraw[black] (12,0) circle (2pt);
\filldraw[black] (16,0) circle (2pt);

\node(u1) at (17,1.8) {$u_{1}$};
 \node(u2) at (19,1.3) {$u_{2}$};

\draw[thick] (16,0) -- (17,1.5)--(19,1)--(16,0);
\filldraw[black] (17,1.5) circle (2pt);
\filldraw[black] (19,1) circle (2pt);

\node(u4) at (17,-1.9) {$u_{4}$};
 \node(u5) at (19,-1.4) {$u_{5}$};

\draw[thick] (16,0) -- (17,-1.5)--(19,-1)--(16,0);
\filldraw[black] (17,-1.5) circle (2pt);
\filldraw[black] (19,-1) circle (2pt);

\node(u21) at (13,-3.3) {$u_{21}$};
 \node(v3) at (11,-1.8) {$v_{3}$};

\draw[thick] (12,0) -- (13,-3)--(11,-1.5)--(12,0);
\filldraw[black] (13,-3) circle (2pt);
\filldraw[black] (11,-1.5) circle (2pt);

 \node(u31) at (16,-2.8) {$u_{31}$};

\draw[thick] (16,0) -- (16,-2.5);
\filldraw[black] (16,-2.5) circle (2pt);
\node(a) at (3.5,-6) {$a)\,G'^*$};
\node(b) at (14.5,-6) {$b)\,G$};
\end{tikzpicture}
\caption{Hypergraph $G'^*$ obtained from $C_3(2,0,0)$ by $v_3$ shrinking on $e_3$ and graph $G$ obtained by replacing each edge of $G'^*$ with a complete graph of corresponding order.}
\label{fig:2}
\end{figure}   
    Let $V_1=\{v_3,u_{21},u_{22},\cdots,u_{2(k-2)}\}$. Since $v_2$ is a cut vertex of $G$, any walk $W\in W_s(G;v_3)$ can be decomposed into two distinct closed walks $W_1$, a $(v_3,v_3)$-walk whose internal vertices are chosen only from $V_1\cup\{v_2\}$ and $W_2$, a $(v_2,v_2)$-walk whose internal vertices are from $V(G)\setminus V_1$. Define a map $f_1$ from $W_s(G;v_3)$ to $W_s(G:u_{11})$ such that $f_1(W)=W_1'W_2$, where $W_1'$ be the walk obtained from $W_1$ by replacing every $v_3$ by $u_{11}$, $u_{21}$ by $v_1$ and $u_{2i},\,i\in[2,k-2]$ by $u_{1i}$, respectively. Then $f_1(W)\in W_s(G;u_{11})$ and $f_1$ is injection, so $M_s(G;v_3)\leq M_s(G;u_{11})$. Also, $M_4(G;v_3)< M_4(G;u_{11})$. Thus $(G;v_3)\prec (G;u_{11}).$\\

    Consider the vertex set $V_2=\{v_1,u_{31},u_{32},u_{33},\ldots,u_{3(k-2)}\}$ and $w\in V_2$. Let $W=W_1W_2$ be a walk in $W_s(G;w,v_3)$, where $W_1$ is the longest $(w,v_2)$-section of the walk $W$ whose internal vertices are chosen from $V(G)\setminus V_1$ and $W_2$ is the remaining $(v_2,v_3)$-section of $W$ whose internal vertices are from $V_1\cup\{v_2\}$. Next, we construct an injective map $f_2$ from $W_s(G;w,v_3)$ to $W_s(G;w,u_{11})$ by $f_2(W)=W'=W_1W_2'$, where $W_2'$ is obtained from $W_2$ by replacing every $v_3$ by $u_{11}$, $u_{21}$ by $v_1$ and $u_{2i},\,i\in[2,k-2]$ by $u_{1i}$, respectively. Hence $(G;w,v_3)\preceq (G;w,u_{11})$.\\

    Let $E_{v_3}=\{wv_3:w\in V_2\}$ and $E_{u_{11}}=\{wu_{11}:w\in V_2\}$. Then by Lemma \ref{EEuv}, $E(G+E_{v_3})< E(G+E_{u_{11}})$. Note that the adjacency matrix of these graphs and the corresponding hypergraph are the same, so the Estrada index. Therefore, from Theorem \ref{secmaxEEG} and Lemma \ref{lemma 4.2}, we get     $EE(C_3(\dfrac{n}{k-1}-3,0,0))<EE(C_2(\dfrac{n}{k-1}-3,1))$.  \end{proof}

\begin{figure}[H]
    \centering
\begin{tikzpicture}[scale=1.15]
\draw[thick] (-0.2,0) ellipse (2.7 and 0.7);
\node at (-0.2,1) {$e_1$};
\filldraw[black] (-2.2,0) circle (2pt);
\node at (-2.2,-0.3) {$u_0$};
\filldraw[black] (-1.45,0) circle (2pt);
\node at (-1.45,-0.3) {$u_{1,1}$};
\filldraw[black] (-0.7,0) circle (2pt);
\node at (-0.7,-0.3) {$u_{1,2}$};
\filldraw[black] (-0.4,0) circle (1pt);
\filldraw[black] (-0.2,0) circle (1pt);
\filldraw[black] (0,0) circle (1pt);
\filldraw[black] (0.35,0) circle (2pt);
\node at (0.35,-0.3) {$u_{1,k-3}$};
\filldraw[black] (1.25,0) circle (2pt);
\node at (1.25,-0.3) {$u_{1,k-2}$};
\filldraw[black] (2,0) circle (2pt);
\node at (2,-0.3) {$u_1$};
\filldraw[black] (2.75,0) circle (2pt);
\node at (2.75,-0.3) {$u_{2,1}$};
\filldraw[black] (3.4,0) circle (2pt);
\node at (3.4,-0.3) {$u_{2,2}$};
\filldraw[black] (3.65,0) circle (1pt);
\filldraw[black] (3.85,0) circle (1pt);
\filldraw[black] (4,0) circle (1pt);
\filldraw[black] (4.3,0) circle (2pt);
\node at (4.3,-0.3) {$u_{2,k-3}$};
\filldraw[black] (5.25,0) circle (2pt);
\node at (5.25,-0.3) {$u_{2,k-2}$};
\filldraw[black] (6,0) circle (2pt);
\node at (6,-0.3) {$u_{2}$};
\draw[thick] (4,0) ellipse (2.5 and 0.7);
\node at (4,1) {$e_2$};
\filldraw[black] (6.75,0) circle (2pt);
\node at (6.75,-0.3) {$u_{3,1}$};
\filldraw[black] (7.5,0) circle (2pt);
\node at (7.5,-0.3) {$u_{3,2}$};
\filldraw[black] (7.8,0) circle (1pt);
\filldraw[black] (8,0) circle (1pt);
\filldraw[black] (8.2,0) circle (1pt);
\filldraw[black] (8.5,0) circle (2pt);
\node at (8.5,-0.3) {$u_{3,k-3}$};
\filldraw[black] (9.45,0) circle (2pt);
\node at (9.45,-0.3) {$u_{3,k-2}$};
\filldraw[black] (10.2,0) circle (2pt);
\node at (10.2,-0.3) {$u_{3}$};
\draw[thick] (8.2,0) ellipse (2.7 and 0.7);
\node at (8.2,1) {$e_3$};
\end{tikzpicture}
\caption{Path $P_3^k$}
\label{fig.4}
\end{figure}
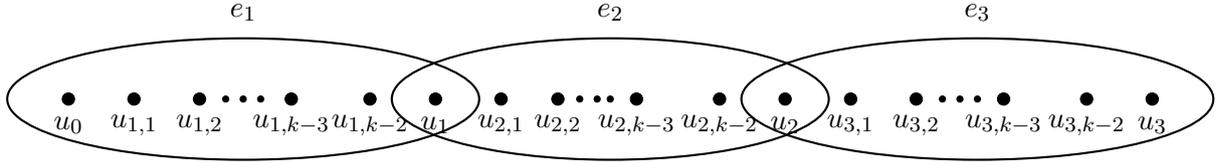

From Theorem \ref{secmaxEEG}, Lemmas \ref{geq3},\ref{C3ni},\ref{lemma 4.2} and \ref{lemma4.3}, we have the following result.
\begin{theorem}\label{unicyclicsecmax}
    Among the $k$-uniform unicyclic hypergraphs of order $n>3(k-1)$, the hypergraph $C_2(\dfrac{n}{k-1}-3,1)$ has the second largest Estrada index.
\end{theorem}
The following corollary is an immediate consequence of Theorems \ref{unicyclicmax} and \ref{unicyclicsecmax}.
\begin{corollary}
    Let $G^*$ be the $k$-uniform hypergraph with maximum Estrada index in $U(n)$ with diameter $d$. Then $G^*\cong C_2(\frac{n}{k-1}-2,0)$ if $d=2$, and $G^*\cong C_2(\frac{n}{k-1}-3,1)$ if $d=3$.
\end{corollary}
\begin{remark}\label{remark4.11}

Let $C_{3,k}$ and $G^{**}$ be $k$-uniform hypergraphs with $n=3(k-1)$ vertices as illustrated in Figure~\ref{fig.3}. Consider the path $P_3^k$ (see Figure~\ref{fig.4}). Let $e_0 = e_1 \setminus \{u_0\}$. 

Since $u_2$ is a cut vertex of $P_3^k$, any closed walk $W\in W_s(P_3^k; u_3)$ can be decomposed as $W = W_1 W_2$, where $W_1$ is a $(u_3,u_3)$-section and $W_2$ is the longest possible $(u_2,u_2)$-section of $W$. Define a map $f_1$ from $ W_s(P_3^k; u_3)$ to $W_s(P_3^k; u_{2,1})$ such that $f_1(W)=W_1'W_2$, where $W_1'$ is obtained from $W_1$ by replacing every occurrence of $u_3$ with $u_{2,1}$, $u_{3,1}$ with $u_1$, and each $u_{3,i}$ with $u_{2,i}$, $i \in [2,k-2]$. Thus, 
\(
M_s(P_3^k, u_3) \leq M_s(P_3^k, u_{2,1}).
\)
Also,
\(
M_4(P_3^k, u_3) < M_4(P_3^k, u_{2,1}),
\)
which implies that $(P_3^k; u_3) \prec (P_3^k; u_{2,1})$.

Now, let $u_{1,1}$ be a pendant vertex of $e_0$. Consider a walk $W = W_1W_2 \in W_s(P_3^k; u_{1,1}, u_3)$ such that $W_1$ is the longest $(u_{1,1}, u_2)$-section and $W_2$ is the remaining $(u_2, u_3)$-section of the walk $W$. We define an injective map $f_2$ from$W_s(P_3^k; u_{1,1}, u_3)$ to $ W_s(P_3^k; u_{1,1}, u_{2,1})$ by $f_2(W)=W_1W_2'$, where $W_2'$ is constructed from $W_2$ by replacing $u_3$ with $u_{2,1}$, $u_{3,1}$ with $u_1$, and $u_{3,i}$ with $u_{2,i}$ for $i \in [2,k-2]$ of $W_2$. This implies:
\(
(P_3^k; u_{1,1}, u_3) \preceq (P_3^k; u_{1,1}, u_{2,1}).
\)

Hence, by Lemma~\ref{EEuv}, we conclude
\[
EE(C_{3,k}) = EE\left(P_3^k - e_1 + (e_0 \cup \{u_3\})\right) < EE\left(P_3^k - e_1 + (e_0 \cup \{u_{1,1}\})\right) = EE(G^{**}).
\]

\end{remark}
\section{Conclusion}

The spectral properties and Estrada index of \( k \)-uniform hypergraphs are examined in this paper. We established bounds for the sum of the \( t \)-largest eigenvalues of a \( k \)-uniform hypergraph and characterized those $k$-uniform hypergraphs that possess exactly two distinct eigenvalues, thereby contributing to the understanding of their spectral structure. Furthermore, we derived bounds for the Estrada index of \( k \)-uniform hypergraphs in terms of order, size and energy. Moreover, we identified the hypergraph that attains the second largest Estrada index within the class of $k$-uniform unicyclic hypergraphs.

\section*{Declarations}
On behalf of all authors, the corresponding author states that there is no conflict of interest.

\bibliographystyle{plain}
\bibliography{ref}
\end{document}